\documentclass[a4paper,12pt,fleqn]{article}

\usepackage{graphicx}      
\usepackage{natbib}        

\usepackage{url}
\usepackage{amsfonts,latexsym}
\usepackage{amssymb}
\usepackage{amsmath}
\usepackage{color}

\usepackage{algorithmic}
\usepackage{algorithm}

\newtheorem{theorem}{Theorem}

\newtheorem{definition}{Definition}

\newtheorem{approximation}{Approximation}
\newtheorem{remark}{Remark}

\newcommand{\msx}{\left. x\right|_1^N}
\newcommand{\msw}{\left. w\right|_1^N}

\newcommand{\D}{\mathrm{d}}
\newcommand{\Prob}{\mathrm{Pr}}
\newcommand{\vol}{\mathrm{vol}}

\newcommand{\ve}{\varepsilon}
\newcommand{\beq}{\begin{equation}}
\newcommand{\eeq}{\end{equation}}
\newcommand{\bea}{\begin{eqnarray}}
\newcommand{\eea}{\end{eqnarray}}
\newcommand{\beas}{\begin{eqnarray*}}
\newcommand{\eeas}{\end{eqnarray*}}
\newcommand{\ba}{\begin{array}}
\newcommand{\ea}{\end{array}}
\newcommand{\bit}{\begin{itemize}}
\newcommand{\eit}{\end{itemize}}
\newcommand{\ben}{\begin{enumerate}}
\newcommand{\een}{\end{enumerate}}
\newcommand{\Viol}{\mathrm{Viol}}

\newcommand{\dss}{\displaystyle}
\newcommand{\Real}[1]{ { {\mathbb R}^{#1} } }

\newcommand{\tran}{ ^{\top} }
\newcommand{\Bp}{\B_{p}}
\newcommand{\Xp}{\mathcal{X_{+}}}
\newcommand{\Xt}{\mathcal{A}}
\newcommand{\xp}{x_{+}}

\renewcommand{\S}{\mathcal{S}}
\newcommand{\X}{\mathcal{X}}
\newcommand{\Y}{\mathcal{Y}}

\newcommand{\E}{\mathcal{E}}

\newcommand{\B}{\mathcal{B}}
\newcommand{\W}{\mathcal{W}}
\newcommand{\V}{\mathcal{V}}
\renewcommand{\P}{\mathcal{P}}
\newcommand{\U}{\mathcal{U}}

\newcommand{\Pd}{\mathbb{P}_\sigma}
\newcommand{\unif}[1]{ {{\lambda}}_{#1} }

\newcommand{\A}{\mathcal{A}}

\newcommand{\ANAS}{{\A}_{\rm NAS}}
\newcommand{\APAS}{{\A}_{\rm PAS}}

\usepackage[textheight=235mm,textwidth=160mm]{geometry} 
\title{\bf Randomized Approximations of the Image Set of Nonlinear Mappings with Applications to Filtering}

\begin{document}


%
%

\author{Fabrizio Dabbene$^1$ Didier Henrion$^{2,3,4}$\\ Constantino Lagoa$^5$ Pavel Shcherbakov$^6$}

\footnotetext[1]{CNR-IEIIT, Politecnico di Torino, Italy}
\footnotetext[2]{CNRS, LAAS,  Toulouse, France}
\footnotetext[3]{Univ. de Toulouse,  Toulouse, France}
\footnotetext[4]{Fac. of Electrical Engr., Czech Tech. Univ.  Prague,  Czech Rep.}
\footnotetext[5]{Electrical Engr. Dept. The Pennsylvania State University,  USA}
\footnotetext[6]{Inst. of Control Sciences, RAS, Moscow, Russia}

\date{Version of \today}

\maketitle

\begin{abstract}
The aim of this paper is twofold: In the first part, we leverage recent results on scenario design to develop randomized algorithms
for approximating the \textit{image set} of a nonlinear mapping, that is, a (possibly noisy) mapping of a set via a nonlinear function.
We introduce minimum-volume approximations which have the characteristic of guaranteeing a low probability of violation, i.e.,
we admit for a probability that some points in the image set are not contained in the approximating set,
but this probability is kept below a pre-specified threshold~$\ve$.
In the second part of the paper, this idea is then exploited to develop a new family of randomized prediction-corrector filters.
These filters represent a natural extension and rapprochement of Gaussian and set-valued filters,
and bear similarities with modern tools such as particle filters.\\[1em]
{\bf Keywords:} Randomized algorithms, filtering, nonlinear systems, semialgebraic sets.\\[1em]
This research was partly funded by CNR-CNRS bilateral project No.~134562.
\end{abstract}


\section{Introduction}
\label{sec-intro}

In recent years, randomized algorithms have gained increasing popularity in the field of control of uncertain systems;
e.g., see~\cite{CaDaTe:11,TeCaDa:13}, due to their ability of dealing with large  and complex uncertainty structures,
thus extending the applicability of the robust control methods. This is obtained by shifting from the robustness paradigm,
where one looks for \textit{guaranteed} performance which should hold for every possible instance of the uncertainty,
to an approach where \textit{probabilistic guarantees} are accepted, i.e.,\ performance is guaranteed only within a given level
of probability $\ve>0$.
The main technical tool that permits to obtain computationally tractable solutions are randomized algorithms,
which could be seen as extensions of the classical Monte Carlo method, and are based on the extraction of random samples of the uncertainty.

\vspace{-.04in}
In this paper, we exploit these ideas for finding reliable approximations of the image of a set through a nonlinear mapping,
which we refer to as the \textit{image set}.
This set is in general nonconvex (possibly not connected), so that classical approximations may be rather difficult to compute
and in general may turn out to be very conservative.

\vspace{-.04in}
The first part of the paper adapts and significantly extends recent results presented by some of the authors in \cite{DaLaSh:10},
where a new definition of ``goodness'' of approximation was provided in probabilistic terms. Namely, an approximating set~$\A$
of a set $\X$ is deemed to be ``good'' if it contains ``most'' of the points in~$\X$ with high probability.
Contrary to classical approximating sets, which are generally either \textit{outer} or \textit{inner},
the ensuing approximating set ``optimally describes'' the set without neither containing nor being contained in it.
This new concept allows to obtain generally tighter approximations, providing a probabilistic  characterization of the set, which is particularly appealing in many contexts (even if it may not be desired in others, such e.g. as safety analysis).
We recall that outer bounding sets, that is, approximations which are guaranteed to contain the set~$\X$, have been very popular
in the set-membership approach, and have been used in designing set-theoretic state estimators for uncertain
discrete-time nonlinear dynamic systems; e.g., see\ \cite{ElGCal:01,ABRC:08}.
Inner ellipsoidal approximations have been introduced for instance in the context of nonlinear programming problems~\citep{NesNem:94}
and tolerance design \citep{WojVla:93}.

\vspace{-.04in}
In Sections II and III, the results of~\cite{DaLaSh:10} are particularized to the specific problem at hand, and are also generalized
by considering a new family of approximating sets which are based on the construction of minimum volume \textit{polynomial superlevel sets},
recently introduced in~\cite{DabHen:13}. With these approximating sets, the original convexity requirement can be relaxed,
since they can be nonconvex/nonconnected, thus allowing for better descriptions.

\vspace{-.04in}
The second part of the paper extends these ideas to the design of probabilistic filters for nonlinear discrete time
dynamical systems subject to random uncertainty. In particular, the image set approximation is used to design a probabilistically optimal predictor filter.
The rationale behind this approach is to ``describe'' the state position at step~$k$ by a set where the state lies with high probability.
The prediction step is then combined with a correction step where the propagated set is trimmed based on the available measurements
at time~$k$ and on the measurement noise assumptions.
A detailed discussion on the proposed randomized prediction correction filter is given in Section IV,
where a discussion on the related literature is also reported.
Numerical examples conclude the paper.


\textbf{Notation:}
{
For a symmetric matrix $P$, $P\succ 0$ means that~$P$ is positive definite.
We denote by $\Bp$ the unit ball in the $\ell_{p}$ norm:
$
    \Bp \doteq \left\{z\in\Real{n}\, :\,\|z\|_p\le 1
        \right\}.
$
The volume or, more precisely, the Lebesgue measure of a compact set $\X$ is denoted by
$\vol\,\X \doteq \int_{\X} \D x.
$
The uniform measure over $\X$ is denoted by $\unif{\X}$, i.e.\ $\unif{\X} $
is such that, for any set $\Y\subseteq \X$,  $\unif{\X}(\Y)=\vol\, \Y /\vol\, \X$.
The set of all polynomials of order less than or equal to $\sigma$ is denoted by $\Pd$.
The monomial basis for this set is represented by the (column) vector $\pi_\sigma(x)$ and any polynomial $q \in  \Pd$ can be expressed in the form
$q(x) = \pi_\sigma\tran(x) q = \pi\tran_{\lceil \sigma/2 \rceil}(x) Q \pi_{\lceil \sigma/2 \rceil}(x)
$ where $q$ is a  vector and $Q$ is a symmetric matrix of appropriate dimensions, referred to as the Gram matrix.
}

\section{Minimum volume approximations}
\vspace{-.1in}
Consider the following nonlinear mapping:
\beq
\label{eq:plus}
\xp  =  f\bigl( x,w \bigr)
\eeq
which represents the one-step evolution of a discrete-time dynamical system,
with $x,\xp\in\Real{n}$ representing the current and future states, respectively,
and $w\in\Real{n_w}$ describing a process noise vector.
We assume that the current state is confined within a compact set $\X\subset\Real{n}$ and that the noise~$w$
also belongs to a given compact set~$\W$.
The problem we are interested in is to find a good approximation to the set of points that can be obtained from~\eqref{eq:plus} by
starting from $x\in\X$ and accounting for all possible values of the noise $w\in\W$, that is, find an approximation of the \textit{image set} defined as
\[
\Xp \!\doteq\! f\bigl(\!\X,\W\!\bigr) \!=\! \bigl\{\xp\!\!\in\Real{n}\!\colon\, \exists x \!\in\! \X, w \!\in\! \W\!\colon \xp \!\!=\!  f( x,w) \bigr\}.
\]

\begin{figure}[h!]
\centerline{\includegraphics[width=0.5\textwidth]{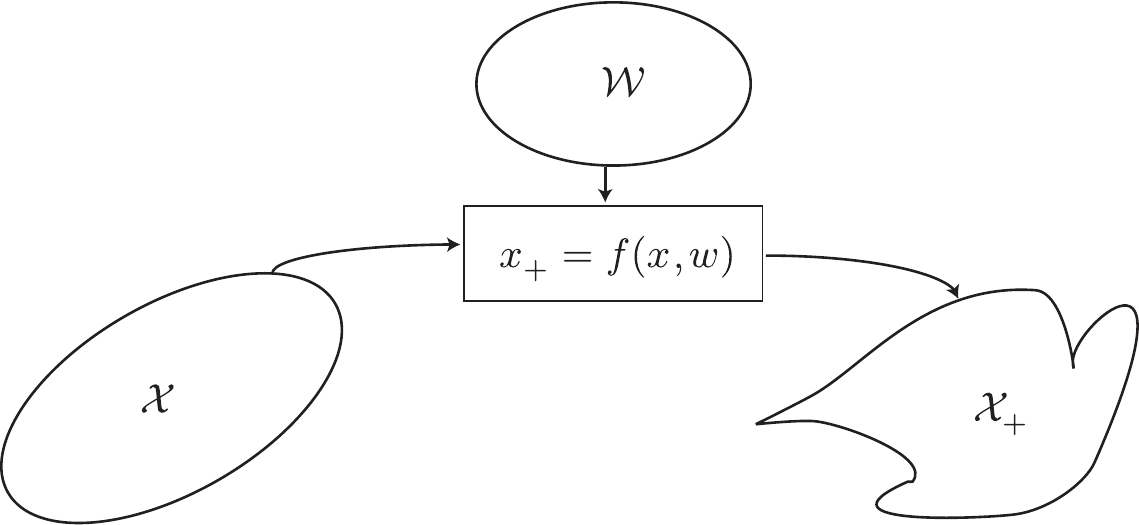}}
\caption{Set image of  $\X$ with noise \!$w\in\W$ under map \!\eqref{eq:plus} .}
\label{fig:Xplus}
\end{figure}
%
%

\vspace{-.03in}
Motivated by the computational complexity issues and the fact that deterministic formulations of the
approximation problem may not be suitable in practical situations,
in~\cite{DaLaSh:10} the authors proposed an original approach for tackling the problem, based on a probabilistic viewpoint.
To this end, a probabilistic information over~$\X$ and~$\W$ is assumed to be known (which is actually the case in many applications),
and an approximation~$\Xt$ of~$\Xp$ is deemed to be ``good''
if it contains all points in~$\Xp$ with high probability.
More formally, we assume that the sets~$\X$ and~$\W$ are endowed
with probability measures~$\mu_{\X}$ and~$\mu_{\W}$. Then, for a given \textit{reliability level}~$\ve$,  the following concept of $\ve$-probabilistic approximation is introduced.

\begin{definition}[$\ve$-probabilistic approximation of $\Xp$]
The

\vspace{-.07in}
set $\Xt$ is an \textit{$\ve$-probabilistic approximation of the set $\Xp$} if
$$
\Viol(\Xt) \le \ve,
$$
\vspace{-.1in}
with
\vspace{-.06in}
\bea
\Viol(\Xt) &\doteq&\Prob\bigl\{x\in\X,w\in\W\colon\, \xp=  f\bigl( x,w \bigr)\not\in \Xt\bigr\}\nonumber\\
&=&\int_{\Xp\setminus\Xt}
f(x,w) \D\mu_{\X}(x) \D\mu_{\W}(w). \label{eq:probviol}
\eea
\end{definition}

\vspace{-.04in}
Note that the probability in \eqref{eq:probviol} is measured with respect to the underlying measures $\mu_{\X}$ and $\mu_{\W}$.
The left-hand side of the equation is referred to as the \textit{violation probability of the set}~$\Xt$.

\vspace{-.04in}
The main characteristic of this approach is that an $\ve$-probabilistic approximating set has neither to cover nor to be fully contained in~$\X$;
it just has to guarantee that the violation probability of the set~$\Xt$ is bounded by~$\ve$.
Clearly, we are interested in finding the smallest among such sets. In the sequel, we first define the two families of approximating sets
considered in this paper.

\vspace{-.04in}
\subsection{Convex approximating sets: Ellipsoids, parallelotopes, and hyperrectangles}

\vspace{-.1in}
The following general description of the \textit{norm-based approximating set} (NAS) was
introduced in \cite{DaLaSh:10}:
\beq
\label{eq:NAS}
   \A(c,P) \doteq \left\{x\!\in\!\Real{n}\colon\, \|P(x-c)\|_p\le \!1,\, P\!=\!P\tran\!\succeq 0
        \right\}.
\eeq
Note that the family of sets above is parameterized by the positive-definite \textit{shape matrix} $P$,
and by the center $c\in\Real{n}$; it represents a generalization of the classical ellipsoidal set for norms different
from the Euclidean one.

\vspace{-.03in}
Indeed, for $p=2$, we obtain the ellipsoid
\[
\E(c,P) \doteq \left\{x\in\Real{n}\colon \, x = c + P^{-1} z, \,  \|z\|_2\le 1
        \right\},
\]
and for $p=\infty$ we get a so-called \textit{elementary parallelotope}, a special-type polytope
with parallel faces (these can be viewed as a particular class of zonotopes \citep{AlBrCa:05} with positive definite generator matrix). In particular, if~$P$ is chosen to be diagonal, we obtain a classical hyper-rectangle.

\vspace{-.03in}
It follows that, in general, the problem of finding the \textit{minimum volume} NAS containing $\Xp$ can be rewritten in the form of the following
robust convex problem:
\beq
\label{NAS}
\begin{array}{l}
    \displaystyle{\min_{P,c}} \log\det(P^{-1}) \\
~~~~~\mbox{s.t.}~~~P=P\tran\succ 0, \\
~~~~~~~~~~~~\|P\xp-c\|_p\le 1  \; \forall \xp\in\Xp.
\end{array}
\tag{NAS-robust}
\eeq

\vspace{-.03in}
\subsection{Nonconvex approximations: Polynomial superlevel sets}
\vspace{-.1in}
In \cite{DabHen:13}, nonconvex set approximations based on the superlevel set of a multidimensional polynomial have been introduced,
and shown to represent a simple and efficient way for describing complex shaped sets.
Formally, assume we are given a compact semialgebraic set
\[
\S:=\{x \in {\mathbb R}^n : b_i(x) \geq 0, \: i=1,2,\ldots,m_b \}
\]
such that $\Xp\!\subseteq is$, with $b_i(x)$ being given polynomials (the set $\S$ is usually a hyper-rectangle). Then, given a polynomial of degree $\sigma\!>0$, i.e.,\ $q \in \Pd$, its polynomial superlevel set $\U(q)$ is defined as
\[
\U(q) \doteq \{ x \in\S : q(x) \geq 1 \}.
\]
In the sequel, we will refer to this family of sets as \textit{polynomial-based approximating sets} (PAS).
The goal is then to find a polynomial $q \in \Pd$ whose corresponding PAS $\U(q)$ contains $\Xp$ and has minimal volume.
In \cite{DabHen:13}, the following problem was formulated:
\begin{equation}
\label{prob:ellone}
\begin{array}{rcll}
&& \displaystyle{\min_{q \in \Pd}}\!\!\! &\! \|q\|_1 \\
&& \mathrm{s.t.} & q \geq 0 \:\:\mathrm{on}\:\: \S ,\\
&&& q \geq 1 \:\:\mathrm{on}\:\: \Xp,
\end{array}
\end{equation}
and the $L^1$ measure above was shown to be a good surrogate to the volume of the set $\U(q)$.
Moreover, theoretically rigorous convergence results were provided in~\cite{DabHen:13}, showing that, if $\Xp$ is semialgebraic,
a hierarchy of outer approximations obtained by solving the above problem for
increasing values of $\sigma$ converges in volume (or, equivalently, almost uniformly) to the set $\Xp$.

\vspace{-.03in}
We now make the important observation that this result still holds for generic compact sets,
not necessarily of semialgebraic form.
Hence, we propose the following family of PASs, where we further relax
the positivity requirement of $q$ over $\S$ by requiring that the polynomial can be expressed as
a \textit{sum-of-squares} (SOS) over~$\S$:
\begin{equation}
\label{PAS-rob}
\begin{array}{rcll}
&& \displaystyle{\min_{q \in \Pd}} \!\!\!&\! \int_\S q(x)\, d x\\
&& \mathrm{s.t.} & q(\xp) \text{~~~is SOS over } \S, \\
&&& q(\xp) \geq 1\quad\forall\xp\in\Xp.
\tag{PAS-robust}
\end{array}
\end{equation}
Note that the requirement that $q$ is SOS  over $\S$ can be formulated using Putinar's Positivstellensatz and
a hierarchy of finite-dimensional convex LMI relaxations which are
linear in the coefficients of $q$. More specifically, we write
$q= r_0 + \sum_i^m r_i b_i$, where $r_0, r_1, \ldots, r_m$ are
polynomial sum-of-squares of given degree $\sigma_{r}$, to be found.
For simplicity, we choose $\sigma_{r}=\sigma$.
For each fixed degree, the problem of finding such polynomials
is an LMI; e.g., see \cite[Section 3.2]{Lasserre:11}.
Note again that problem \eqref{PAS-rob} is a robust convex program (that is, for every fixed $\xp$, we get a convex program).

\vspace{-.04in}
\section{Randomized approximations}

\vspace{-.07in}
In this section, we propose a simple randomized algorithm to construct, with arbitrarily high probability,
$\ve$-probabilistic approximations of the image set~$\Xp$.
The suggested randomized procedure is rather straightforward: We draw~$N$ samples of $x\in\X$
according to the measure $\mu_{\X}$, obtaining the multisample
\beq
\label{eq:multis_x}
\msx\doteq\{x^{(1)},\dots x^{(N)}\}
\eeq
and, similarly, we draw $N$ instances of the random noise $w\in\W$, drawn according to the (given)
noise measure $\mu_{\W}$
\beq
\label{eq:multis_w}
\msw\doteq\{w^{(1)},\dots w^{(N)}\}.
\eeq
Based on these samples, we compute the points
\[
\xp^{(i)}=f(x^{(i)},w^{(i)}), \quad i=1,\ldots,N.
\]
Then, we construct the approximating set $\Xt$ as the minimum-size bounding set containing
the random points $\left. \xp\right|_1^N$ thus obtained.
This is doable by solving one of the following two scenario problems\footnote{%
Note the slight abuse of notation adopted: $\ANAS$ denotes the
NAS set defined by the optimal values of $P$ and $c$ obtained by solving problem \eqref{NAS-scen}. Same for $\APAS$.}:
\vspace{-.04in}
\beq
\label{NAS-scen}
\begin{array}{l}
\ANAS\,\doteq\,\arg\displaystyle{\min_{P,c}} \dss\log\det(P^{-1}) \\
  ~~~~\mbox{s.t.}~~P=P\tran\succ 0, \\
 ~~~~~~~~~\|P\xp^{(i)}-c\|_p\le 1,   \;   i=1,\ldots,N
\tag{NAS-random}
\end{array}
\eeq
\vspace{-.05in}
and
\vspace{-.05in}
\begin{equation}
\label{PAS-scen}
\begin{array}{l}
\APAS\,\doteq\,\arg\displaystyle{\min_{q \in \Pd}}  \dss \int_\S q(x)\, d x\\
 ~~~~\mathrm{s.t.}~~~q(\xp) \text{~~~is SOS over } \S, \\
 ~~~~~~~~~~~q(\xp^{(i)}) \geq 1, \;        i=1,\ldots,N.
\tag{PAS-random}
\end{array}
\end{equation}

\begin{figure}[h!]
\centerline{\includegraphics[width=0.5\textwidth]{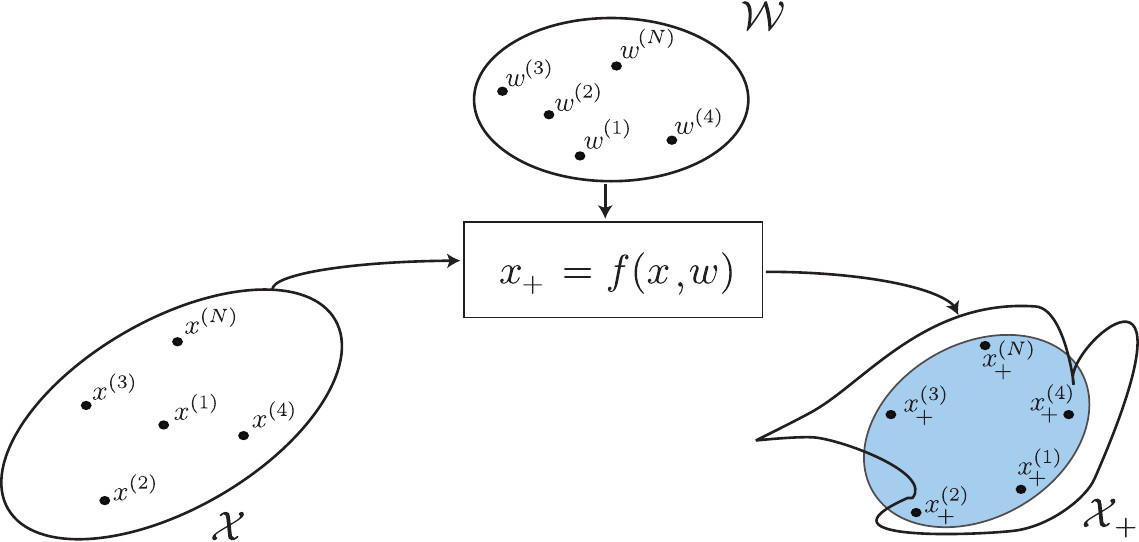}}
\caption{Randomized approximation of the image set.}
\label{fig:Xplus-scen}
\end{figure}

\vspace{-.07in}
The idea at the basis of the approach is depicted in Fig.~\ref{fig:Xplus-scen}, where random samples of~$x$ and~$w$ are drawn
and mapped via the nonlinear function $f(\cdot,\cdot)$.
Then, an approximation of $\Xp$ is obtained by constructing a minimum volume approximating set around these points.

\vspace{-.12in}
Note that both problems above represent simple convex optimization problems which can be efficiently solved by available SDP solvers.
In particular: (a)~problem \eqref{NAS-scen} reduces to a linear program (LP) for $p=1$ and diagonal $P$;
(b)~the constraints  $q(\xp^{(i)}) \!\geq\! 1$ in problem \eqref{PAS-scen} can also be expressed as an LP in the coefficients of $q(x)$.
The procedure is summarized in the following algorithm.

\vspace{-.04in}
\begin{algorithm}[h!]
\caption{\!Randomized Image Set \! Approximation
\label{RSIA}}
\begin{algorithmic}[1]
\REQUIRE $\X$,
\ENSURE $\ANAS$ (or $\APAS$)
\STATE draw $x^{(1)},\ldots,x^{(N)}$ over $\X$ and
$w^{(1)},\ldots,w^{(N)}$ over $\W$
\STATE construct $\xp^{(i)}=f(x^{(i)},w^{(i)}), \quad i=1,\ldots,N$
\STATE solve problem \eqref{NAS-scen} (or \eqref{PAS-scen})
\STATE return $\ANAS$ (or $\APAS$).
\end{algorithmic}
\end{algorithm}

\vspace{-.07in}
The key features of the algorithm above are as follows: (i)~there is no need in knowing the probability measures
$\mu_{\X}$ and $\mu_{\W}$ explicitly, only the capability of obtaining random samples from these are mandatory for running the algorithm;
(ii)~formal results
on the probabilistic quality of the obtained approximations are available.
In particular, the properties of the algorithm follow from the well-known results of \cite{CalCam:06tac}
on the so-called scenario optimization, which  allow to bound \textit{a priori} (before the algorithm is run)
the violation probability of the solution.

\begin{theorem}
Let $\A$ (i.e., $\ANAS$ or $\APAS$) be the outcome of Algorithm~1.
Fix a \textit{risk level} $\ve>0$ and let~$d$ be the ``size'' of the optimization problem \ref{NAS-scen} (or \ref{PAS-scen}).
Then we have
\vspace{-.04in}
\beq
\label{onestep-prob}
\Pr\bigl\{\!
\Viol({\A})\!>\!\ve\!
\bigr\}\le
\Phi(\ve,N)\doteq
\sum_{j=0}^{d-1}\!
{\!N\choose{j}\!} \ve^j(1\!-\!\ve)^{N-j}.
\eeq
Moreover, the above bound is tight, see \cite{CamGar:08,Calafiore:10}.
\end{theorem}

\vspace{-.04in}
Note that the result above makes explicit use of the size~$d$ of the convex problem.
This represents the number of free variables in the convex optimization problem; see Table~\ref{table-d}.

\vspace{-.05in}
\begin{table}[h!]
\caption{Number of design variables\label{table-d}}
\begin{center}
\begin{tabular}{|c|c|}
\hline
\hline
Approximating set & Number of design variables \\ \hline
\hline
NAS - Ellipsoid/Parallelotope  & $ d=\frac{n(n+1)}{2}+n$\\ \hline
NAS - Hyperrectangle  & $d=2n$\\ \hline
PAS - polynomial of degree $\sigma$ & $d={n_{x} \choose {n_{x}+\sigma}}$\\
\hline\hline
\end{tabular}
\end{center}
\label{default}
\end{table}%

\vspace{-.03in}
In order to obtain an expression for the desired sample size~$N$ (the number of samples ensuring that the designed approximation
is $\ve$-optimal with probability at least $\delta$),  we need to invert equation~\eqref{onestep-prob}.
The best bound currently available is proposed in \cite{ATLR:15}, where it is shown that, given $\ve,\delta\in(0,1)$,
if $N$ is chosen as
\[
N\ge  \frac{e}{e-1}\frac{1}{\ve} \left(d+\ln \frac{1}{\delta}\right),
\]
then $\A$ is an $\ve$-probabilistic approximation of $\Xp$ with probability at least $(1-\delta)$.

\vspace{-.04in}
\subsection{Numerical Example}

\begin{figure}[h!]
\centerline{\includegraphics[width=0.5\textwidth]{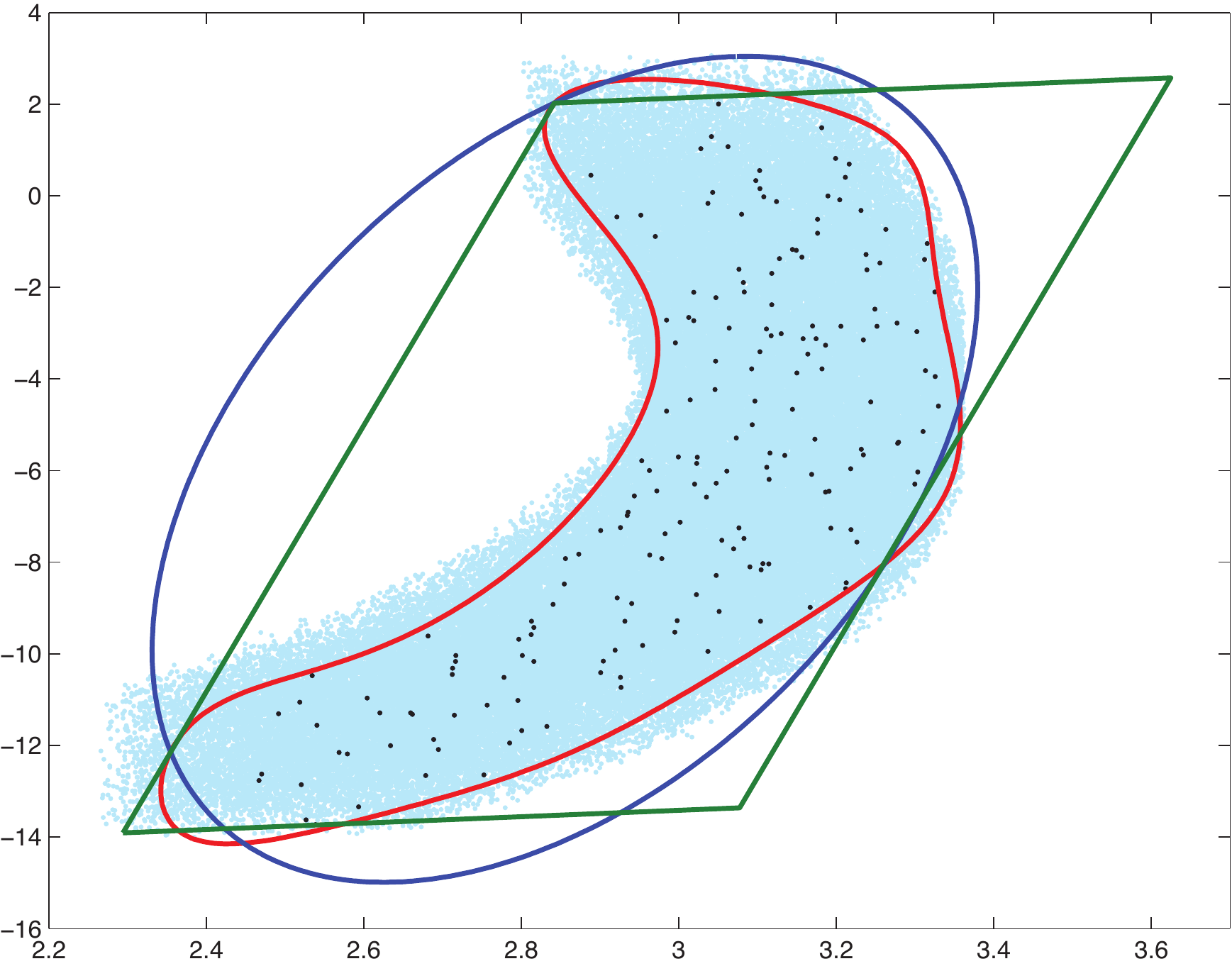}}
\caption{Image set of  \eqref{sysF}, with ellipsoidal NAS (blue), parallelotopic NAS (green), and PAS (red) approximations.}
\label{fig:image}
\end{figure}

\vspace{-.07in}
To show the approximating sets obtained by Algorithm 1, we consider the nonlinear mapping
\vspace{-.05in}
\bea
\label{sysF}
\nonumber
\xp(1) & = &
\sin x(2)+3\cos x(2)+w(1),\\
\xp(2) & = &  3x(1)-20\log\bigl(1+x(2)\bigr)+w(2),
\eea
where $x\sim\unif{\X}$ with $\X=[0,1]^{2}$ and  $w\sim\unif{\W}$ with $\W=[-0.2,0.2]^{2}$.
Figure \ref{fig:image} depicts $N=200$ randomly generated points in~$\Xp$, and the corresponding NAS and PAS probabilistic approximations.
The shaded area represents the ``true'' set~$\Xp$ obtained by generating one million points in it.
These approximations are indeed seen to be neither inner nor outer ones.

\vspace{-.04in}
\section{Randomized Filtering}

\vspace{-.07in}
The results of the previous sections are here reinterpreted in a filtering setting.
The idea is rather intuitive: One can apply the image set description in an iterative way
and use the approximation as a description of the state position at time $k$.
For the sake of simplicity, we limit our exposition to the case where NAS descriptions are employed.

\vspace{-.04in}
Consider the following discrete-time (autonomous) nonlinear system
\vspace{-.05in}
\bea
\label{system-x}
x_{k + 1} & = & f\bigl( x_{k},w_{k} \bigr),\\
\label{system-meas}
y_{k}     & = & g\bigl( x_{k}\bigr) +v_{k},
\eea
\noindent
Equation (\ref{system-x}) represents the systems dynamics, and $x_{k}$, $x_{k+1}\in\Real{n}$ denote the state at time instants $k$ and $k+1$,
respectively, while $w_{k}\in\Real{n_w}$ is a noise vector, usually referred to as process noise. 
Equation (\ref{system-meas}) represents a nonlinear measurement mapping associated with the system equations, where $y_{k}\in\Real{n_y}$ is the measurement vector and $v_k\in\Real{n_v}$ is the measurement noise vector.

\vspace{-.04in}
Various approaches have been adopted in the literature to address the filtering problem for this type of systems.
In the Extended Kalman filtering approach, and, in general, in all the subsequently developed Gaussian Filtering techniques,
(e.g., see\ \cite{ItoXio:00}), probabilistic information on the states and noise is considered.
The main approximation made in this setup is that \textit{all information about the state and noise is collected in the first two moments of this distribution};
in other words, the state and (update and measurement) noise can be described by a Gaussian distribution.
Note that, even if at $k=0$ everything could indeed be Gaussian, this is clearly not true after the first step,
since Gaussian distributions are generally not preserved under nonlinear mappings.
Approaches of this line are also those based on Unscented Kalman Filtering,\; e.g., see~\cite{Sarkka:07}.

\vspace{-.04in}
A philosophically different approach is the one based on set-valued filtering, in which no probabilistic assumptions are made on the
nature of the noise and on the \emph{a priori} information on $x$; instead, a so-called unknown-but-bounded approach is used.
In those works (e.g., see\ \cite{Calafiore:05}), the \emph{a priori} information is given in terms of a bounded set
where the initial state is guaranteed to lie, and at each time step both process and measurement noise
are assumed to take values from given bounded sets.
Then, all information on the predicted state is itself captured by a set (usually an ellipsoid) which is
\textit{guaranteed to contain the state at time $k+1$}.

\vspace{-.04in}
In this work, we propose a rapprochement between those two models and assume that the process and measurement noise are random,
with given distributions over bounded support sets.
That is, we assume that, at each time instant $k$, $w_k \in \W$ and $v_k \in \V$, where $\W,\V$ are bounded sets%
\footnote{We may assume different $\W_{k},\V_{k}$ for different time instants.}.
Moreover, we assume that $w_{k}$ has stochastic nature, and we denote by $\P_{\W}$
the probability measure on $\W$.
Note that no stochasticity assumption is made on $v_{k}$.

\vspace{-.04in}
The next approximation is crucial in our setup, and, from the philosophical point of view, it plays the same role
of the Gaussian approximation commonly made in Gaussian filtering, in which the first two moments of the distribution
are taken as representatives of the whole pdf of the propagated state.

\begin{approximation}[State description]
At time $k$, \textit{all} state information is captured by a NAS set $\A_{k}=\A(P_{k},c_{k})$
(e.g.,\ an ellipsoid) with center $c_{k}$ and shape matrix $P_{k}$.
Moreover, we assume that $x_{k}$ has stochastic nature, and it is \textit{uniformly distributed} over $\A_{k}$.
\end{approximation}

\vspace{-.06in}
Note again  that this constitutes one of the main approximations introduced in our setup:
At each step we `reset' the density information about the state and make the implicit assumption
that $x_k$ has a \textit{uniform distribution} $\lambda$ with ellipsoidal support $\A_{k}$.

\vspace{-.02in}
\subsection{Probabilistically Guaranteed  Simulation}

\vspace{-.07in}
In this section, we concentrate on the state dynamics equation only and
propose a prediction filter for simulating the behavior of system (\ref{system-x}).
To this end, assume that at step~$k$, the state~$x_k$ is guaranteed (in probability)
to lie inside a NAS set, that is
$
x_k\in\A_k\doteq\A(x_k,P_k).
$
Then, we can apply the set image approximation algorithm; namely,
given a probabilistic parameter~$\ve$, we choose $N$ so as to guarantee a desired approximation level,
and we start by drawing uniform random samples
$x_k^{(1)},\ldots,x_k^{(N)}\sim \unif{\A_k}$, and $N$ noise samples
$w_k^{(1)},\ldots,w_k^{(N)}\sim \unif{\W}$.
Based on these samples, we compute the points
$x_{k+1}^{(i)}=f(x_k^{(i)},w_k^{(i)})$, $i=1,\ldots,N$, and solve the \ref{NAS-scen} problem.
The proposed procedure is summarized in Algorithm~\ref{RA-pred}.
\begin{algorithm}[h!]
\caption{One-step Randomized Prediction Filter
\label{RA-pred}}
\begin{algorithmic}[1]
\REQUIRE $\A_k$,
\ENSURE $\A_{k+1}$
\STATE draw  $x_k^{(1)},\ldots,x_k^{(N)}\sim\unif{\A_k}$ and $w_k^{(1)},\ldots,w_k^{(N)}\sim \unif{\W}$
\STATE construct $x_{k+1}^{(i)}=f(x_k^{(i)},w_k^{(i)}), \quad i=1,\ldots,N$
\STATE solve problem \ref{NAS-scen} 
\STATE return $\A_{k+1}$.
\end{algorithmic}
\end{algorithm}

\vspace{-.05in}
\begin{remark}[$m$-step ahead prediction]
The iterative one-step prediction scheme discussed above has been
introduced having in mind the prediction-correction scheme that we formulate in the next section.
However, there are cases when it could be convenient to directly perform an $m$-step ahead prediction,
in one shot. To do this, a simple solution is to ``propagate'' the state samples $m$-times ahead
in the following way: Generate $N$ samples $x_k^{(i)}$, $i=1,\ldots,N$, and $mN$ samples
of uncertainty $w_{k+j-1}^{(i)}$, $i=1,\ldots,N$, $j=1,\ldots,m$, and iteratively compute
\[
x_{k+j}^{(i)}=f(x_{k+j-1}^{(i)},w_{k+j-1}^{(i)}), \quad i=1,\ldots,N;\, j=1,\ldots,m.
\]
Then, we can directly construct a NAS approximation of $\A_{k+m}$ based on the multisample
$x_{k+m}^{(1)},\ldots,x_{k+m}^{(N)}$.
\end{remark}

\vspace{-.06in}
\subsection{Measurement update}

\vspace{-.07in}
If a measurement is available, it can be immediately used to \textit{reject} incorrect predictions,
as discussed in Algorithm~\ref{RPCF} (RPCF) below.

\begin{algorithm}[h!]
\caption{One-step Randomized Prediction-Correction Filter\! (RPCF)
\label{RPCF}}
\begin{algorithmic}[1]
\REQUIRE $\A_k$,
\ENSURE $x_{k+1}$
\STATE draw  $x_k^{(1)},\ldots,x_k^{(N)}\sim\unif{\A_k}$
and $w_k^{(1)},\ldots,w_k^{(N)}\sim\unif{\W}$
\STATE construct $x_{k+1}^{(i)}=f(x_k^{(i)},w_k^{(i)}), \quad i=1,\ldots,N$
\STATE construct $z_{k+1}^{(i)}=y_{k+1}- g(x_{k+1}^{(i)}), \quad i=1,\ldots,N$
\IF {$z_{k+1}^{(i)}\not\in\V_{k}$}
\STATE reject $x_{k+1}^{(i)}$
\ENDIF
\STATE solve problem \ref{NAS-scen} 
\STATE return $\A_{k+1}$.
\end{algorithmic}
\end{algorithm}

\begin{remark}[Similarities with particle filtering]
In the so called \textit{particle filtering}, particles are generated and fed through the dynamical system.
Then, they are kept/dis\-carded based on Bayesian-type reasonings. Our setup has strong similarities.
Indeed, one can construct a ``particle-like'' implementation of the predictor-corrector randomized algorithm,
where at each step not all samples are generated again, but instead, the non-discarded ones are fed through.
This algorithm seems to provide a good practical behavior. However, it still needs a thorough theoretical analysis,
which will be the subject of future work.
\end{remark}

\vspace{-.06in}
\section{Numerical Example}

\vspace{-.08in}
As a filtering example, we tested the specific two-dimensional nonlinear system with scalar measurements,
borrowed from \cite{ABRC:08}, which was considered in the context of set filtering by means of zonotopic approximation sets:
\vspace{-.08in}
\bea
\label{system_ex}
 \label{ex_one}
x_{k + 1}(1) & = & -0.7x_{k}(2) + 0.1x_k(2)^2+ 0.1x_k(1)x_k(2)+ 0.1 {\rm e}^{x_k(1)} + w_k(1),\\  \label{ex_two}
x_{k + 1}(2) & = & x_k(1) + x_k(2) - 0.1x_k(1)^2 +0.2x_k(1)x_k(2) + w_k(2), \\ \label{y_meas}
y_{k}       & = & x_k(1) + x_k(2) + v_{k}, 
\eea
where the propagation noise is $|w_k(1)|\leq 0.1$, $|w_k(2)|\leq 0.1$;
the measurement error is $|v_k|\leq 0.2$, $k>0$; the initial state belongs to the box $\X_{0} = [-3,\, 3]\times[-3,\,3]$.
Note that the measurements are linear in both $x$ and $v$, and the dynamics is linear in $w$.

\vspace{-.02in}
To shape the experiment, we first generated a ``realized'' trajectory $x^{sim}_{k}$ and the respective sequence of
measurements $y_{k}$, $k=1,\dots,K$. Then, given $\A_k$ at step~$k$, we sample~$N$ points in $\A_k$ and propagate
them via \eqref{ex_one}--\eqref{ex_two} to obtain a set~$\tilde\X_k$ of  candidate next-step points $x^{(i)}_{k+1}$.
We call $x^{(i)}_{k+1}\in \tilde\X_k$  {\it good} if it fits the measurement, i.e.
\beq
\label{slab}
|y_{k+1} - x_{k+1}^{(i)}(1) - x_{k+2}^{(i)}(2)| \leq  0.2,
\eeq
see (\ref{y_meas}). We call this account for the measurements as {\it correction} step or {\it rejection}
(of bad points). Because of the linearity of the measurements, the good points are those who fall in the
{\it slab}~$\S_{k+1}$ defined by~(\ref{slab}).

\vspace{-.04in}
Next, we build a minimum volume set $\A_{k+1}$ (of the shape adopted) that contains only the good points and adopt it as a
probabilistic approximation set for~$x_{k+1}$.

\begin{figure}[h!]
\centerline{\includegraphics[width=0.5\textwidth]{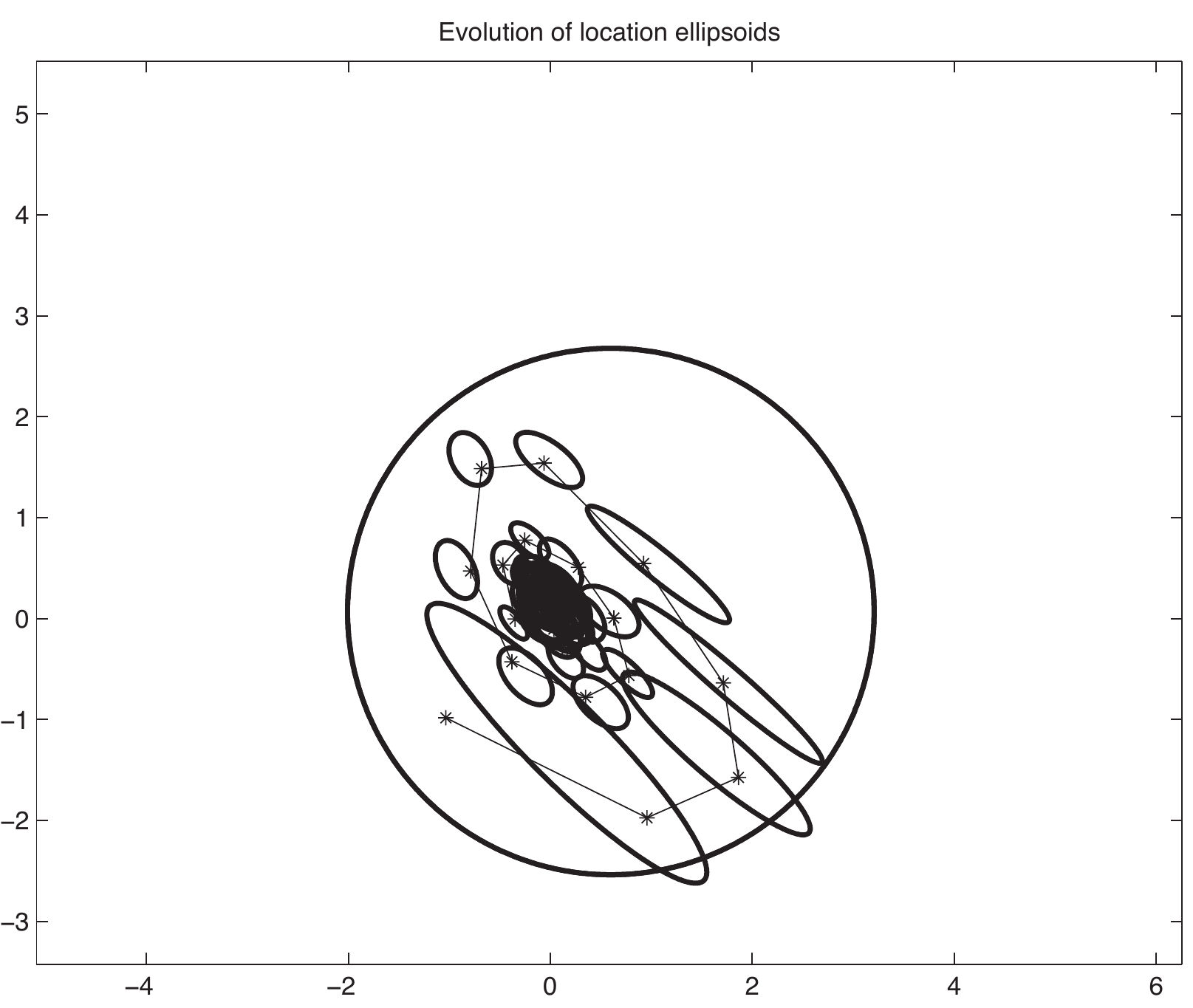}}
\caption{\small Evolution of the approximation ellipses.}
\label{fig:locell_1}
\end{figure}

\vspace{-.04in}
Figure~\ref{fig:locell_1} represents the evolution of the sets $\A_{k}$ in the case of ellipsoidal shapes (the circle represents the
ellipse~$\E([0.6,0.07]\tran, 6.8I)$ of initial uncertainty). The stars correspond to the reference trajectory $x^{sim}_{k}$.

%
%

\begin{figure}[h!]
\centerline{\includegraphics[width=0.45\textwidth]{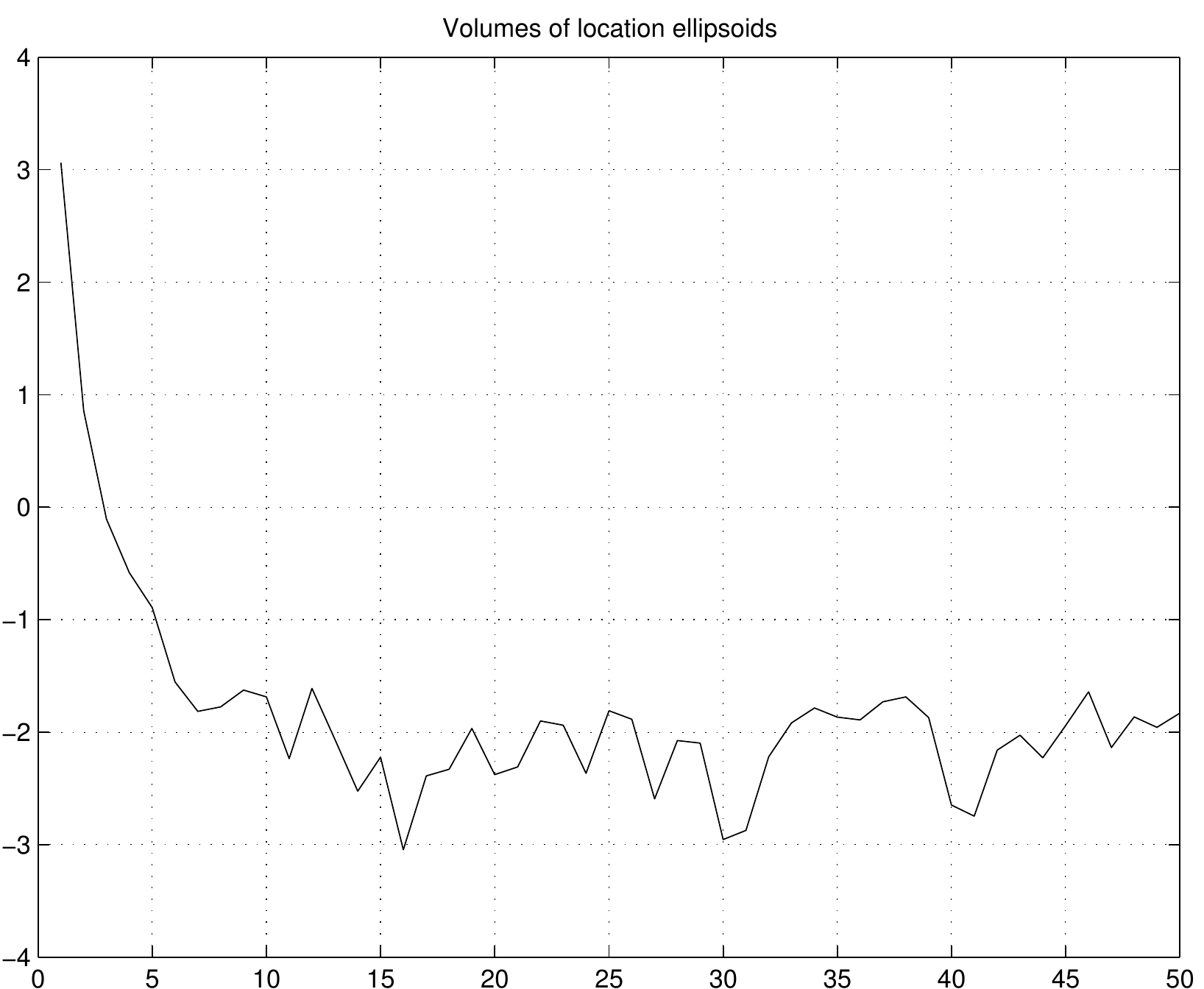}~\includegraphics[width=0.45\textwidth]{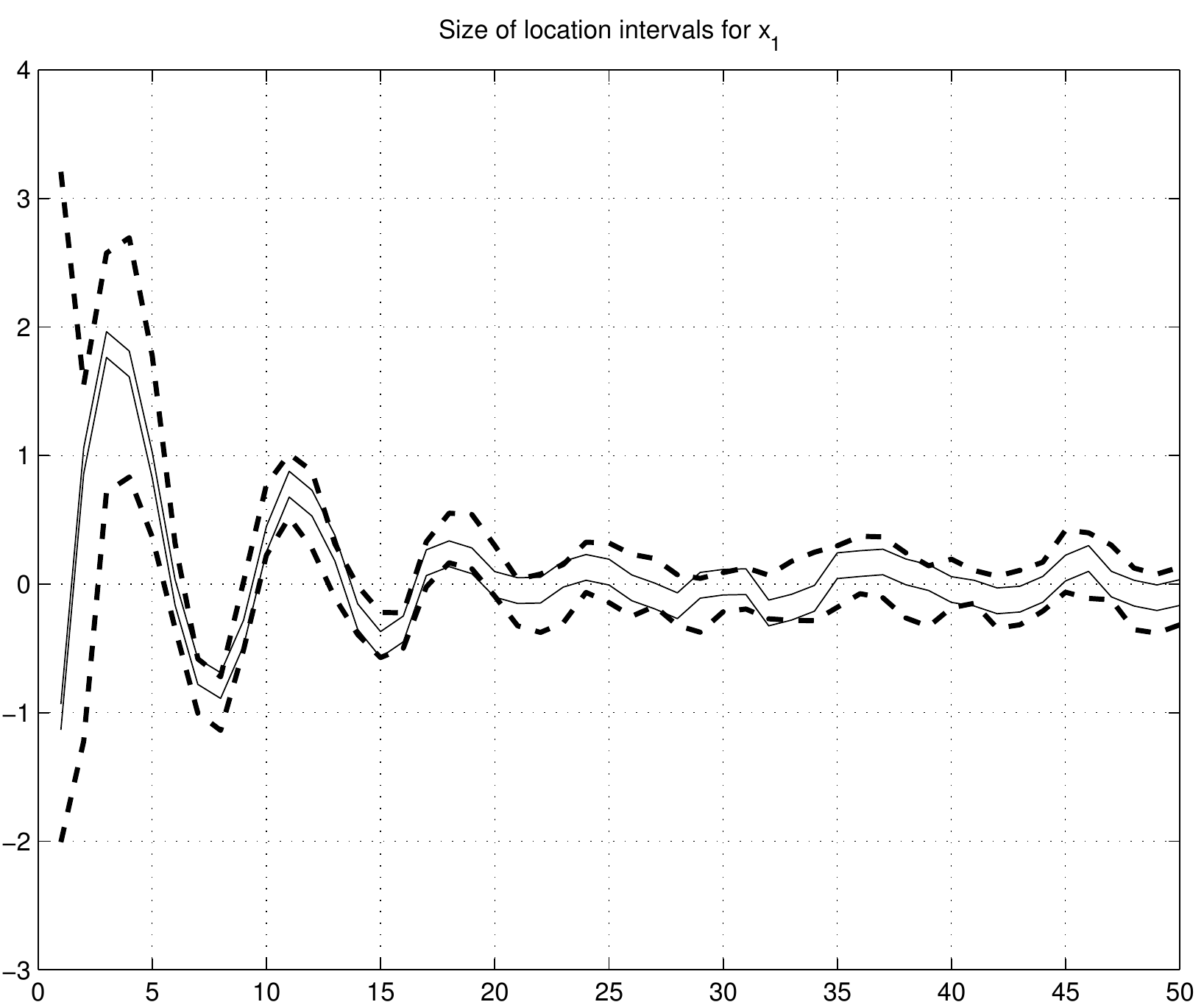}}
\caption{\small The volumes of the approximation ellipsoids (left) and estimation error for $x^{sim}_1(k)$ (right).}
\label{fig:vollell_1}
\end{figure}

For the same realization we plot $\log\bigl(\vol(\A_k)\bigr)$, see Fig.~\ref{fig:vollell_1} (left).
We remark that this quantity is of the same order of magnitude (or better) than that for the zonotopes in~\cite{ABRC:08}, but of course the computations are
way simpler, and this advantage will get stronger as the dimension of~$x$ grows.

%

The filtering goal was to estimate $x_k(1)$, the first component of $x_{k}$.
In Fig.~\ref{fig:vollell_1} (right), we
plot the $x(1)$-span of our ellipses, representing the estimation error, (bold dashed curves),
while with the solid line we plot the worst-case bounds on $x(1)$ obtained from the reference trajectory:
$x^{sim}_k(1)\pm 0.1$ (see definition (\ref{ex_one})--(\ref{ex_two}) keeping in mind that the bounded propagation noise~$|w_k(1,2)|\leq 0.1$
enters additively in the dynamics).

\vspace{-.04in}
Notice that sometimes the inner corridor goes outside the bounds obtained by our method. This is exactly what
illustrates the nature of our approach: We do neglect low-probability events.

\begin{remark}[Comments on re-sampling and re-use]
In the RPCF Algorithm, many samples could be rejected due to the measurement equation. Experiments seem to show that in such cases it
is preferable to \textit{re-sample} the samples rejected in the measurement phase,
so as to guarantee that the new set $\A_{k+1}$ is constructed based on a sufficiently large number of samples.
Another point worth noticing is that, on the contrary, having built the $\A_k$ at step $k$, we already possess $N_{good}\leq N$ points inside it,
so that they can be re-used (propagated), and just the rest $N-N_{good}$ ones are to be newly sampled; this is called {\it re-use}.
\end{remark}

\vspace{-.03in}
\section{Conclusions}
\vspace{-.05in}
We developed an efficient algorithm for constructing $\ve$-probabilistic approximations of the image set of a noisy nonlinear mapping.
The recursive application of such technique, combined with an ad-hoc rejection procedure, leads to the design of a new family of randomized
prediction-correction filters, which showed themselves pretty competitive with both classical set-theoretical and Gaussian filtering methods.
The derivations of the RCPF were limited to the NAS setting, but we stress that the use of PAS descriptions for the construction
of prediction/correction filters is also possible in principle.
However, we admit that in that case, the computations may become too cumbersome for online implementations.


\begin{thebibliography}{17}
\providecommand{\natexlab}[1]{#1}
\providecommand{\url}[1]{\texttt{#1}}
\expandafter\ifx\csname urlstyle\endcsname\relax
  \providecommand{\doi}[1]{doi: #1}\else
  \providecommand{\doi}{doi: \begingroup \urlstyle{rm}\Url}\fi

\bibitem[Alamo et~al.(2005)Alamo, Bravo, and Camacho]{AlBrCa:05}
T.~Alamo, J.M. Bravo, and E.F. Camacho.
\newblock Guaranteed state estimation by zonotopes.
\newblock \emph{Automatica}, 41\penalty0 (6):\penalty0 1035--1043, 2005.

\bibitem[Alamo et~al.(2008)Alamo, Bravo, Redondo, and Camacho]{ABRC:08}
T.~Alamo, J.M. Bravo, M.J. Redondo, and E.F. Camacho.
\newblock A set-membership state estimation algorithm based on {DC}
  programming.
\newblock \emph{Automatica}, 44\penalty0 (1):\penalty0 216--224, 2008.

\bibitem[Alamo et~al.(2015)Alamo, Tempo, Luque, and Ramirez]{ATLR:15}
T.~Alamo, R.~Tempo, A.~Luque, and {D.R.} Ramirez.
\newblock Randomized methods for design of uncertain systems: Sample complexity
  and sequential algorithms.
\newblock \emph{Automatica}, 52:\penalty0 160--172, 2015.

\bibitem[Calafiore and Campi(2006)]{CalCam:06tac}
G.~Calafiore and M.~Campi.
\newblock The scenario approach to robust control design.
\newblock \emph{IEEE Transactions on Automatic Control}, 51\penalty0
  (5):\penalty0 742--753, 2006.

\bibitem[Calafiore et~al.(2011)Calafiore, Dabbene, and Tempo]{CaDaTe:11}
G.~Calafiore, F.~Dabbene, and R.~Tempo.
\newblock Research on probabilistic methods for control system design.
\newblock \emph{Automatica}, 47:\penalty0 1279--1293, 2011.

\bibitem[Calafiore(2005)]{Calafiore:05}
G.C. Calafiore.
\newblock Reliable localization using set-valued nonlinear filters.
\newblock \emph{IEEE Transactions on Systems, Man, and Cybernetics},
  35\penalty0 (2):\penalty0 189--197, 2005.

\bibitem[Calafiore(2010)]{Calafiore:10}
G.C. Calafiore.
\newblock Random convex programs.
\newblock \emph{{SIAM} Journal on Optimization}, 20\penalty0 (6):\penalty0
  3427--3464, 2010.

\bibitem[Campi and Garatti(2008)]{CamGar:08}
M.~Campi and S.~Garatti.
\newblock The exact feasibility of randomized solutions of robust convex
  programs.
\newblock \emph{{SIAM} Journal on Optimization}, 19:\penalty0 1211---1230,
  2008.

\bibitem[Dabbene and Henrion(2013)]{DabHen:13}
F.~Dabbene and D.~Henrion.
\newblock Set approximation via minimum-volume polynomial sublevel sets.
\newblock In \emph{Proc. of the European Control Conference}, 2013.

\bibitem[Dabbene et~al.(2010)Dabbene, Lagoa, and Shcherbakov]{DaLaSh:10}
F.~Dabbene, C.~Lagoa, and P.~Shcherbakov.
\newblock On the complexity of randomized approximations of nonconvex sets.
\newblock In \emph{Proc. IEEE Multiconference on Systems and Control}, 2010.

\bibitem[{El~Ghaoui} and Calafiore(2001)]{ElGCal:01}
L.~{El~Ghaoui} and G.~Calafiore.
\newblock Robust filtering for discrete-time systems with bounded noise and
  parametric uncertainty.
\newblock \emph{IEEE Transactions on Automatic Control}, 46\penalty0
  (7):\penalty0 1084--1089, July 2001.

\bibitem[Ito. and Xiong(2000)]{ItoXio:00}
K.~Ito. and K.~Xiong.
\newblock Gaussian filters for nonlinear filtering problems.
\newblock \emph{IEEE Transactions on Automatic Control}, 45\penalty0
  (5):\penalty0 910--927, 2000.

\bibitem[Lasserre(2015)]{Lasserre:11}
J.B. Lasserre.
\newblock Level sets and non {G}aussian integrals of positively homogeneous
  functions.
\newblock \emph{Int. Game Theory Review}, 17\penalty0 (1), 2015.

\bibitem[Nesterov and Nemirovski(1994)]{NesNem:94}
Y.~Nesterov and A.S. Nemirovski.
\newblock \emph{Interior Point Polynomial Algorithms in Convex Programming}.
\newblock SIAM, Philadelphia, 1994.

\bibitem[Sarkka(2007)]{Sarkka:07}
S.~Sarkka.
\newblock On unscented {K}alman filtering for state estimation of
  continuous-time nonlinear systems.
\newblock \emph{IEEE Trans. on Automatic Control}, 52\penalty0 (9):\penalty0
  1631--1641, 2007.

\bibitem[Tempo et~al.(2013)Tempo, Calafiore, and Dabbene]{TeCaDa:13}
R.~Tempo, G.C. Calafiore, and F.~Dabbene.
\newblock \emph{Randomized Algorithms for Analysis and Control of Uncertain
  Systems: With Applications}.
\newblock Springer, 2nd edition, 2013.

\bibitem[Wojciechowski and Vlach(1993)]{WojVla:93}
J.~M. Wojciechowski and J.~Vlach.
\newblock Ellipsoidal method for design centering and yield estimation.
\newblock \emph{Transactions on Computer-Aided Design of Integrated Circuits
  and Systems}, 12:\penalty0 1570--1579, 1993.

\end{thebibliography}
\end{document}